\documentclass[11pt]{amsart}%
\usepackage{amssymb,amsmath,amsfonts,latexsym,amsthm,geometry,graphicx}
\usepackage[normalem]{ulem}
\usepackage{amsmath}
\usepackage{amsfonts}
\usepackage{color, soul}
\usepackage{amssymb}
\usepackage{tikz-cd}
\usepackage[all]{xy}
\usepackage{graphicx}%
\providecommand{\U}[1]{\protect\rule{.1in}{.1in}}
\newtheorem{theorem}{Theorem}[section]
\newtheorem{corollary}[theorem]{Corollary}

\theoremstyle{definition}
\newtheorem{problem}[theorem]{Problem}

\newtheorem{definition}[theorem]{Definition}
\begin{document}
\title{Nonlinear variants of a theorem of Kwapie\'{n}}
\author[Renato Macedo]{Renato Macedo}
\address{Departamento de Matem\'{a}tica \\
Universidade Federal da Para\'{\i}ba \\
58.051-900 - Jo\~{a}o Pessoa, Brazil.}
\email{renato\_burity@hotmail.com}
\author[Daniel Pellegrino]{Daniel Pellegrino}
\address{Departamento de Matem\'{a}tica \\
Universidade Federal da Para\'{\i}ba \\
58.051-900 - Jo\~{a}o Pessoa, Brazil.}
\email{dmpellegrino@gmail.com}
\author[Joedson Santos]{Joedson Santos}
\address{Departamento de Matem\'{a}tica \\
Universidade Federal da Para\'{\i}ba \\
58.051-900 - Jo\~{a}o Pessoa, Brazil.}
\email{joedsonmat@gmail.com or joedson@mat.ufpb.br}
\thanks{2010 Mathematics Subject Classification: 46B28, 46T99, 47B10, 54E40}
\thanks{Renato Macedo is partially supported by Capes, Daniel Pellegrino is supported
by CNPq Grant 307327/2017-5 and Grant 2019/0014 Para\'{\i}ba State Research
Foundation (FAPESQ) and Joedson Santos is supported by CNPq Grant
309466/2018-0 and Grant 2019/0014 Para\'{\i}ba State Research Foundation (FAPESQ)}
\keywords{Lipschitz operators, almost summing operators, nonlinear operators}
\maketitle

\begin{abstract}
A famous result of S. Kwapie\'{n} asserts that a linear operator from a Banach
space to a Hilbert space is absolutely $1$-summing whenever its adjoint is
absolutely $q$-summing for some $1\leq q<\infty$; this result was recently
extended to Lipschitz operators by Chen and Zheng. In the present paper we
show that Kwapie\'{n}'s and Chen--Zheng theorems hold in a very relaxed
nonlinear environment, under weaker hypotheses. Even when restricted
to the original linear case, our result generalizes Kwapie\'{n}'s theorem
because it holds when the adjoint is just almost summing. In addition, a variant for
$\mathcal{L}_{p}$-spaces, with $p\geq2$, instead of Hilbert spaces is provided.

\end{abstract}

\section{Introduction}

The theory of absolutely summing linear operators was originated from
Grothendieck's Resum\'{e} \cite{gro} and since the publication of the papers
of Lindenstrauss, Pe{\l }czy{\'{n}}ski and Pietsch \cite{LP, Pie}, it has
performed a fundamental role in Banach Space Theory and its applications. In
the last decades the path from the linear to the nonlinear theory of
absolutely summing operators was encouraged by seminal works of Pietsch
\cite{PPPP}, Farmer and Johnson \cite{FJ}, among others. Nowadays, nonlinear
variants of absolutely summing operators have been explored in different
settings with applications in different fields of pure and applied mathematics
(see, for instance, \cite{bayart, perez, vieira} and the references therein).
Let $X,Y$ be Banach spaces over $\mathbb{K}=\mathbb{R}$ or $\mathbb{C}$ and,
from now on, $B_{X}$ denotes the closed unit ball of $X$, and $X^{\ast}$
represents its topological dual. If $1\leq p<\infty,$ a linear operator
$T:X\longrightarrow Y$ is \textit{absolutely} $p$-\textit{summing} if there
exists a constant $C$ such that%
\begin{equation}
\left(
{\displaystyle\sum\limits_{i=1}^{n}}
\left\Vert T(x_{i})\right\Vert ^{p}\right)  ^{1/p}\leq C\Vert(x_{i})_{i=1}%
^{n}\Vert_{p,w},\label{s}%
\end{equation}
for every $x_{1},...,x_{n}\in X$ and all positive integers $n$, where
\[
\Vert(x_{i})_{i=1}^{n}\Vert_{p,w}:=\sup_{\varphi\in B_{X^{\ast}}}\left(
{\displaystyle\sum\limits_{i=1}^{n}}
|\varphi(x_{i})|^{p}\right)  ^{1/p}.
\]
The class of absolutely $p$-summing linear operators from $X$ to $Y$ will be
represented, as it is usual, by $\Pi_{p}\left(  X,Y\right)  $ and the infimum
of all $C$ satisfying (\ref{s}) defines a norm on $\Pi_{p}\left(  X,Y\right)
$, denoted by $\pi_{p}(T)$. The related notion of almost summing operators
will be important for our purposes. According to \cite{BBJ}, a linear operator
$T:X\longrightarrow Y$ is \textit{almost} $p$-\textit{summing} if there is a
constant $C$ such that
\begin{equation}
\left(  {\int\nolimits_{0}^{1}}\left\Vert \sum_{i=1}^{n}r_{i}(t)T(x_{i}%
)\right\Vert ^{2}dt\right)  ^{\frac{1}{2}}\leq C\left\Vert (x_{i})_{i=1}%
^{n}\right\Vert _{p,w},\label{ppk}%
\end{equation}
for every positive integer $n$ and $x_{1},...,x_{n}\in X$. Above, $r_{i}$ are
the Rademacher functions, defined as
\[%
\begin{array}
[c]{cccc}%
r_{i}: & \left[  0,1\right]   & \longrightarrow & \mathbb{R}\\
& t & \longmapsto & \text{sign}\left(  \sin2^{i}\pi t\right)  .
\end{array}
\]
The infimum of the constants $C$ satisfying (\ref{ppk}) is denoted by
$\pi_{al.s.p}(T)$. The class of all almost $p$-summing operators is denoted by
$\Pi_{al.s.p}$. When $p=2$, these operators are simply called \textit{almost
summing} and we write $\Pi_{al.s}$ instead of $\Pi_{al.s.2}$ (see
\cite[Chapter 12]{djt}). By \cite[Proposition 12.5]{djt},
\begin{equation}
\bigcup_{1\leq q<\infty}\Pi_{q}(X,Y)\subseteq\Pi_{al.s}(X,Y).\label{inccc}%
\end{equation}
For details, we refer to the classical book by Diestel, Jarchow and Tonge
\cite{djt}.

From now on $\mathcal{I}$ denotes an operator ideal, $\mathcal{L}$ denotes the
ideal of all bounded linear operators and $u^{\ast}$ represents the adjoint of
a linear operator $u$. We denote, as usual, by $\mathcal{I}^{dual}$ the
\textit{dual ideal} of $\mathcal{I}$ (see \cite[page 186]{djt}), that is,
\[
\mathcal{I}^{dual}(X,Y):=\{u\in\mathcal{L}(X,Y);~u^{\ast}\in\mathcal{I}%
(Y^{\ast},X^{\ast})\}.
\]
In general, the adjoint of an absolutely $p$-summing operator may fail to be
absolutely $p$-summing and \textit{vice-versa}. For instance, considering the
formal identity $i:\ell_{2}\longrightarrow c_{0}$, it follows that $i^{\ast}$
is absolutely $p$-summing for all $1\leq p<\infty,$ whereas $i$ and
$i^{\ast\ast}$ aren't absolutely $p$-summing. Then, the following problem
arises naturally:

\begin{problem}
For which Banach spaces $X$ and $Y$ and $1\leq p,q<\infty$ does the following
inclusion hold
\begin{align*}
\Pi_{q}^{dual}(X,Y) \subseteq\Pi_{p}(X,Y)?
\end{align*}

\end{problem}

A remarkable result due to Kwapie\'{n} \cite{K} shows that a linear operator
$T$ from a Banach space $X$ to a Hilbert space $H$ is absolutely $1$-summing
whenever $T^{\ast}$ is absolutely $q$-summing for some $1\leq q<\infty$. In
other words,

\begin{theorem}
$($see \cite[Theorem 2.21]{djt}$)$\label{kuap} If $X$ is a Banach space and
$H$ is a Hilbert space, then
\begin{align*}
\bigcup_{1\leq q<\infty}\Pi_{q}^{dual}(X,H)\subseteq\Pi_{1}(X,H).
\end{align*}

\end{theorem}

Kwapie\'{n}'s Theorem was recently extended to Lipschitz operators by Chen and
Zheng \cite[Theorem 3.1]{CZ} as we shall see next. Let $X$ be a pointed metric
space with a base point denoted by $0$, and let $Y$ be a Banach space. The
Lipschitz space $Lip_{0}(X;Y)$ is the Banach space of all Lipschitz operators
$T:X\longrightarrow Y$ such that $T(0)=0$, under the Lipschitz norm defined
by
\[
Lip(T)=\sup\left\{  \frac{\Vert T(x)-T(y)\Vert}{d(x,y)}~:~x,y\in X,~x\neq
y\right\}  .
\]
When $Y=\mathbb{K}$, we write $X^{\#}=Lip_{0}(X;\mathbb{K})$. A cornerstone of
the nonlinear theory of absolutely summing operators is the paper of Farmer
and Johnson \cite{FJ}, which introduces the notion of absolutely summing
operators to the Lipschitz framework as follows: a Lipschitz operator $T\colon
X\rightarrow Y$ is \textit{Lipschitz} $p$-\textit{summing i}f there is a
constant $C$ such that
\begin{equation}
\left(  \sum_{i=1}^{n}\Vert T(x_{i})-T(q_{i})\Vert^{p}\right)  ^{\frac{1}{p}%
}\leq C\sup_{\varphi\in B_{X^{\#}}}\left(  \sum_{i=1}^{n}|\varphi
(x_{i})-\varphi(q_{i})|^{p}\right)  ^{\frac{1}{p}} , \label{765765}%
\end{equation}
for all $x_{1},\ldots,x_{n},q_{1},\ldots,q_{n}\in X$ and all positive integers
$n$ (for related papers on nonlinear summing operators we refer to \cite{JA,
JAproc, Fer, MM} and the references therein). According to \cite{Sawashima},
the Lipschitz adjoint $T^{\#}$ of $T\in Lip_{0}(X;Y)$ is defined as the
continuous linear operator $T^{\#}\colon Y^{\#}\longrightarrow X^{\#}$ given
by $T^{\#}(g):=g\circ T$. The variant of Kwapie\'{n}'s Theorem to Lipschitz
operators due to Chen and Zheng reads as follows:

\begin{theorem}
\label{py} $($see \cite[Theorem 3.1]{CZ}$)$ Let $X$ be a pointed metric space
and $H$ be a Hilbert space. If $T\in Lip_{0}(X;H)$ is such that $T^{\#}%
|_{H^{\ast}}$ is absolutely $q$-summing for some $1\leq q<\infty$, then $T$ is
Lipschitz $1$-summing.
\end{theorem}

Following this vein, other questions seem natural:

\begin{problem}
\label{probl3} Are there other nonlinear extensions/generalizations of the
Kwapie\'{n} Theorem?
\end{problem}

\begin{problem}
\label{probl4} Are there variants of the Kwapie\'{n} Theorem in which the
range is not necessarily a Hilbert space?
\end{problem}

In this paper we shall answer these problems. To state our main result we
shall recall the notion of $\mathcal{L}_{p}$-spaces. If $\lambda>1$ and $1\leq
p\leq\infty$, a Banach space $Y$ is an $\mathcal{L}_{p,\lambda}$%
\textit{-space} if every finite dimensional subspace $E$ of $Y$ is contained
in a finite dimensional subspace $Y_{0}$ of $Y$ for which there is an
isomorphism $v:Y_{0}\longrightarrow\ell_{p}^{\dim(Y_{0})}$ such that $\Vert
v\Vert\Vert v^{-1}\Vert\leq\lambda$. When $Y$ is an $\mathcal{L}_{p,\lambda}%
$\textit{-space for a certain }$\lambda$, we just say that $Y$ is an
$\mathcal{L}_{p}$-space$.$ It is simple to observe that Hilbert spaces are
$\mathcal{L}_{2}$-spaces. As a consequence of the main results of the present
paper we conclude that the original result of Kwapie\'{n} can be improved as follows:

\begin{theorem}
\label{nn11}If $X$ is a Banach space and $Y$ is an $\mathcal{L}_{p}$-space and
$2\leq p<\infty$, then%
\[
\Pi_{al.s.p}^{dual}(X,Y)\subseteq\Pi_{1}(X,Y).
\]
In particular, when $Y$ is a Hilbert space%
\begin{equation}
\Pi_{al.s}^{dual}(X,Y)\subseteq\Pi_{1}(X,Y).\label{2222}%
\end{equation}

\end{theorem}

Note that from (\ref{inccc}) it is obvious that (\ref{2222}) recovers the
statement of Kwapie\'{n}'s Theorem. The paper is organized as follows. In
Section 2 we state our main result (Theorem \ref{q1}); in Section 3 we prove
Theorem \ref{q1} and in Section 4 we provide applications of the main result;
for instance, we shall generalize the Chen--Zheng Theorem (Theorem \ref{py})
following the lines of what is done in Theorem \ref{nn11}.

\section{Main result}

We start off by recalling that the sequence of Rademacher functions $\left(
r_{i}\right)  _{i=1}^{\infty}$ is orthonormal; thus
\begin{equation}
\int_{0}^{1}\left\vert \sum_{i=1}^{n}r_{i}(t)a_{i}\right\vert ^{2}%
dt=\sum_{i=1}^{n}|a_{i}|^{2} \label{7654}%
\end{equation}
for all $(a_{i})_{i=1}^{n} \in\ell^{n}_{2}$ and all positive integers $n.$ The
well-known Kahane Inequality shows that the spaces generated by the Rademacher
functions have equivalent $L_{p}$ norms:

\begin{theorem}
[Kahane inequality]Let $0<p,q<\infty$. Then there is a constant $\mathrm{K}%
_{p,q}>0$ for which
\[
\left(  \int_{0}^{1}\left\Vert \sum_{i=1}^{n}r_{i}(t)x_{i}\right\Vert
^{q}dt\right)  ^{\frac{1}{q}}\leq\mathrm{K}_{p,q}\left(  \int_{0}%
^{1}\left\Vert \sum_{i=1}^{n}r_{i}(t)x_{i}\right\Vert ^{p}dt\right)
^{\frac{1}{p}}
\]
holds, regardless of the choice of a Banach space $X$ and of finitely many
vectors $x_{1},\dots,x_{n}\in X$.
\end{theorem}

The previous theorem, combined with (\ref{7654}) recovers the Khinchin
inequality (see \cite[page 10]{djt}), when $q=2$ and $X=\mathbb{K}$, and in
this case $\mathrm{K}_{p,2}^{-1}$ is usually denoted by $A_{p}$; it is simple
to observe that $A_{2}=1$ and it is well-known that $A_{1}=(\sqrt{2})^{-1}$.

In the recent years a series of works (\cite{jfa, joed, adv, LONDON,
PSDianaE}) have shown that several important results of the theory of summing
operators in fact need essentially no linear structure to be valid. The proof
of our main result shall rely on the abstract environment presented in
\cite{adv}.

From now on, unless stated otherwise, $p\in\lbrack1,\infty).$ Let $X$ and $Y$
be non-void sets, $\mathcal{H}\left(  X;Y\right)  $ be a non-void family of
mappings from $X$ to $Y$ and $K$ be a compact Hausdorff space. Let
\[
R\colon K\times X\times X\longrightarrow\lbrack0,\infty)~\text{and}%
\mathrm{~}S\colon\mathcal{H}\left(  X;Y\right)  \times X\times
X\longrightarrow\lbrack0,\infty)
\]
be arbitrary mappings. A mapping $f\in\mathcal{H}\left(  X;Y\right)  $ is
$RS$-\textit{abstract} $p$-\textit{summing} if there is a constant $C$ such
that
\[
\left(  \sum_{i=1}^{n}S(f,x_{i},q_{i})^{p}\right)  ^{\frac{1}{p}}\leq
C\sup_{\varphi\in K}\left(  \sum_{i=1}^{n}R\left(  \varphi,x_{i},q_{i}\right)
^{p}\right)  ^{\frac{1}{p}},
\]
for all $x_{1},\ldots,x_{n},q_{1},\ldots,q_{n}\in X$ and all positive integers
$n$. We define
\[
{\mathcal{H}}_{RS,p}(X;Y)=\left\{  f\in\mathcal{H}\left(  X;Y\right)  :\text{
}f\text{ is }RS\text{-abstract }p\text{-summing}\right\}  .
\]

Let $X$ be a pointed metric space and $Y$ be a Banach space. We shall choose a
suitable Banach space $X^{d}\subseteq\mathbb{K}^{X}$ containing $X^{\#}$ such
that
\[
f^{a}\colon Y^{\#}\longrightarrow X^{d}~,f^{a}(h):=h\circ f
\]
is well-defined for all $f\in\mathcal{H}\left(  X;Y\right)  $. In general, for
our applications, it will be enough to consider $X^{d}=X^{\#}$. Note that,
when $f$ is a linear operator between Banach spaces we have $f^{a}(y^{\ast
})=f^{\ast}(y^{\ast})$ for all $y^{\ast}\in Y^{\ast}$, i.e., $f^{a}$ is a kind of \textit{abstract adjoint} of $f$. We shall need the following
properties of $R$ and $S$:

\medskip

\textbf{(I)} $S(f,x,q)\leq\left\Vert f(x)-f(q)\right\Vert $ for all
$f\in\mathcal{H}\left(  X;Y\right)  $ and $x,q\in X$.

\textbf{(II)} $K\subseteq X^{d}$, $B_{T^{a}(Y^{\ast})}\subseteq K$ and
\[
\left\vert g(x)-g(q)\right\vert \leq R(g,x,q)\ \ \text{for all}\ \ g\in
B_{T^{a}(Y^{\ast})}\text{ and}\ x,q\in X.
\]
In the next section we shall present a nonlinear Kwapie\'{n}-type theorem
valid for $\mathcal{L}_{p}$-spaces instead of just Hilbert spaces. As an
illustration, when the main result of this section is restricted to the linear
case, we conclude that if $X$ is a Banach space and $Y$ is an $\mathcal{L}%
_{p}$-space, with $2\leq p<\infty$, then
\[
\Pi_{al.s.p}^{dual}(X,Y)\subseteq\Pi_{1}(X,Y).
\]
Our techniques are, in general, different from the ones used in the proofs of
Theorems \ref{kuap} and \ref{py}; for instance, unlike what is done in the
aforementioned results, we do not use the Pietsch Domination Theorem. Our main
result is the following:

\begin{theorem}
\label{q1} Let $X$ be a pointed metric space, $Y$ be an $\mathcal{L}_{p}%
$-space, with $2\leq p<\infty$, and $f:X\longrightarrow Y$ an arbitrary map.
If $R$ and $S$ satisfy \textbf{(I)} and \textbf{(II)} and $f^{a}|_{Y^{\ast}}$
is almost $p$-summing, then $f$ is $RS$-abstract $1$-summing.
\end{theorem}

In Section 4, applications of our main result are provided. For
instance, choosing suitable $R,S$ we obtain extensions of the theorems of
Kwapie\'{n} and Chen--Zheng (see Subsections \ref{ss1} and \ref{ss2}).

\section{Proof of the main result}

Let $f\in\Pi_{al.s.p}^{dual}(X,Y)$. Fix $x_{1},...,x_{n},q_{1},...,q_{n}\in X$
and a subspace $Y_{0}$ of $Y$ containing $\{f(x_{1}),...,f(x_{n}%
),f(q_{1}),...,f(q_{n})\}$ and for which there is an isomorphism
$v:Y_{0}\longrightarrow\ell_{p}^{m}$ such that $\Vert v\Vert\Vert v^{-1}%
\Vert\leq\lambda$. Then, using the monotonicity of the $\ell_{p}$ norms, the
Khinchin inequality and denoting by $\left(  e_{k}\right)  _{k=1}^{m}$ the
canonical basis of $\ell_{p^{\ast}}^{m}=\left(  \ell_{p}^{m}\right)  ^{\ast}$,
where $p^{\ast}$ is the conjugate of $p$, it follows that
\begin{align*}
\sum_{i=1}^{n}S(f,x_{i},q_{i})  &  \overset{\mathbf{(I)}}{\leq}\sum_{i=1}%
^{n}\Vert f(x_{i})-f(q_{i})\Vert\\
&  =\sum_{i=1}^{n}\Vert v^{-1}vf(x_{i})-v^{-1}vf(q_{i})\Vert\\
&  \leq\Vert v^{-1}\Vert\cdot\sum_{i=1}^{n}\Vert vf(x_{i})-vf(q_{i})\Vert\\
&  =\Vert v^{-1}\Vert\cdot\sum_{i=1}^{n}\left(  \sum_{k=1}^{m}|e_{k}\left(
vf(x_{i})-vf(q_{i}\right)  )|^{p}\right)  ^{\frac{1}{p}}\\
&  \leq\Vert v^{-1}\Vert\cdot\sum_{i=1}^{n}\left(  \sum_{k=1}^{m}|e_{k}\left(
vf(x_{i})-vf(q_{i})\right)  |^{2}\right)  ^{\frac{1}{2}}\\
&  \leq\sqrt{2}\cdot\Vert v^{-1}\Vert\cdot\sum_{i=1}^{n}{\int\nolimits_{0}%
^{1}}\left\vert \sum_{k=1}^{m}e_{k}\left(  vf(x_{i})-vf(q_{i}\right)  )\cdot
r_{k}(t)\right\vert dt\\
&  =\sqrt{2}\cdot\Vert v^{-1}\Vert\cdot\sum_{i=1}^{n}{\int\nolimits_{0}^{1}%
}\left\vert \sum_{k=1}^{m}e_{k}(vf(x_{i}))\cdot r_{k}(t)-\sum_{k=1}^{m}%
e_{k}(vf(q_{i}))\cdot r_{k}(t)\right\vert dt\\
&  =\sqrt{2}\cdot\Vert v^{-1}\Vert\cdot\sum_{i=1}^{n}{\int\nolimits_{0}^{1}%
}\left\vert \sum_{k=1}^{m}v^{\ast}\left(  e_{k}\right)  (f(x_{i}))\cdot
r_{k}(t)-\sum_{k=1}^{m}v^{\ast}\left(  e_{k}\right)  (f(q_{i}))\cdot
r_{k}(t)\right\vert dt.
\end{align*}
By the Hahn-Banach Theorem we can extend $v^{\ast}(e_{k})$ to $Y^{\ast}$, and
thus
\begin{align*}
\sum_{i=1}^{n}S(f,x_{i},q_{i})  &  \leq\sqrt{2}\cdot\Vert v^{-1}\Vert\cdot
\sum_{i=1}^{n}{\int\nolimits_{0}^{1}}\left\vert \sum_{k=1}^{m}f^{a}v^{\ast
}\left(  e_{k}\right)  (x_{i})\cdot r_{k}(t)-\sum_{k=1}^{m}f^{a}v^{\ast
}\left(  e_{k}\right)  (q_{i})\cdot r_{k}(t)\right\vert dt\\
&  =\sqrt{2}\cdot\Vert v^{-1}\Vert\cdot\sum_{i=1}^{n}{\int\nolimits_{0}^{1}%
}\left\vert \left(  \sum_{k=1}^{m}r_{k}(t)f^{a}v^{\ast}\left(  e_{k}\right)
\right)  (x_{i})-\left(  \sum_{k=1}^{m}r_{k}(t)f^{a}v^{\ast}\left(
e_{k}\right)  \right)  (q_{i})\right\vert dt.
\end{align*}
Combining the previous inequality with \textbf{(II)} we obtain
\[
\sum_{i=1}^{n}S(f,x_{i},q_{i})\leq\sqrt{2}\cdot\Vert v^{-1}\Vert\cdot
{\int\nolimits_{0}^{1}}\left\Vert \sum_{k=1}^{m}r_{k}(t)f^{a}v^{\ast}%
(e_{k})\right\Vert dt\cdot\sup_{g\in K}\sum_{i=1}^{n}R(g,x_{i},q_{i})
\]
and, using the monotonicity of the $L_{p}$ norms, we have
\[
\sum_{i=1}^{n}S(f,x_{i},q_{i})\leq\sqrt{2}\cdot\Vert v^{-1}\Vert\cdot\left(
{\int\nolimits_{0}^{1}}\left\Vert \sum_{k=1}^{m}r_{k}(t)f^{a}v^{\ast}%
(e_{k})\right\Vert ^{2}dt\right)  ^{\frac{1}{2}}\cdot\sup_{g\in K}\sum
_{i=1}^{n}R(g,x_{i},q_{i}).
\]
Since $f^{a}|_{Y^{\ast}}$ is almost $p$-summing, it follows that
$f^{a}|_{Y^{\ast}}v^{\ast}$ also is almost $p$-summing and, since $\Vert
(e_{k})_{k=1}^{n}\Vert_{p,w}=1$ in $\ell_{p^{\ast}}^{m},$ we have
\begin{align*}
\sum_{i=1}^{n}S(f,x_{i},q_{i})  &  \leq\sqrt{2}\cdot\Vert v^{-1}\Vert\cdot
\pi_{al.s.p}(f^{a}|_{Y^{\ast}}v^{\ast})\cdot\sup_{g\in K}\sum_{i=1}%
^{n}R(g,x_{i},q_{i})\\
&  \leq\sqrt{2}\cdot\Vert v^{-1}\Vert\cdot\Vert v^{\ast}\Vert\cdot\pi
_{al.s.p}(f^{a}|_{Y^{\ast}})\cdot\sup_{g\in K}\sum_{i=1}^{n}R(g,x_{i},q_{i})\\
&  =\sqrt{2}\cdot\Vert v^{-1}\Vert\cdot\Vert v\Vert\cdot\pi_{al.s.p}%
(f^{a}|_{Y^{\ast}})\cdot\sup_{g\in K}\sum_{i=1}^{n}R(g,x_{i},q_{i})\\
&  \leq\sqrt{2}\cdot\lambda\cdot\pi_{al.s.p}(f^{a}|_{Y^{\ast}})\cdot\sup_{g\in
K}\sum_{i=1}^{n}R(g,x_{i},q_{i}).
\end{align*}
As a corollary, since every Hilbert space is an
$\mathcal{L}_{2}$-space, we have an abstract generalization of the Kwapie\'{n} Theorem:

\begin{corollary}
[Abstract Kwapie\'{n}-type Theorem]\label{7777}Let $X$ be a pointed metric
space, $H$ be a Hilbert space and $f\in\mathcal{H}(X,H)$. If $R$ and $S$
satisfy \textbf{(I)} and \textbf{(II) }and $f^{a}|_{H^{\ast}}$ is absolutely
$q$-summing for some $1\leq q<\infty$, then $T$ is $RS$-abstract $1$-summing.
\end{corollary}

\section{Applications}

In this section we show that Theorem \ref{q1} provides new Kwapie\'{n}-type
theorems for several classes of linear and nonlinear summing operators (and
recovers the known results).

\subsection{Absolutely summing linear operators\label{ss1}}

Let $X$ and $Y$ be Banach spaces and $T\colon X\longrightarrow Y$ be a bounded
linear operator. Considering the weak-star topology in $B_{X^{\ast}}$, and
letting
\[
\left(  K,\mathcal{H}\left(  X;Y\right)  \right)  =\left(  B_{X^{\ast}%
},\mathcal{L}(X;Y)\right)  ,
\]%
\[
\left\{
\begin{array}
[c]{c}%
R\colon K\times X\times X\longrightarrow\lbrack0,\infty)\\
R(\varphi,x,q)=|\varphi(x-q)|,
\end{array}
\right.
\]
and%
\[
\left\{
\begin{array}
[c]{c}%
S\colon\mathcal{L}\left(  X;Y\right)  \times X\times X\longrightarrow
\lbrack0,\infty)\\
S(T,x,q)=\Vert T(x-q)\Vert,
\end{array}
\right.
\]
it is plain that $T$ is $RS$-abstract $p$-summing if, and only if, there is a
constant $C$ such that
\[
\left(  \sum_{i=1}^{n}\Vert T(x_{i}-q_{i})\Vert^{p}\right)  ^{\frac{1}{p}}\leq
C\sup_{\varphi\in B_{X^{\ast}}}\left(  \sum_{i=1}^{n}|\varphi(x_{i}%
-q_{i})|^{p}\right)  ^{\frac{1}{p}},\label{36M}%
\]
for all $x_{1},\ldots,x_{n},q_{1},\ldots,q_{n}\in X$ and every positive
integer $n$. Thus, $T$ is absolutely $p$-summing if, and only if, it is
$RS$-abstract $p$-summing. Taking $X^{d}=X^{\#}$, it follows that
$T^{a}|_{Y^{\ast}}=T^{\ast}$ and the hypotheses \textbf{(I)} and \textbf{(II)}
are satisfied. If $Y$ is an $\mathcal{L}_{p}$-space, with $2\leq p<\infty$,
Theorem \ref{q1} tells us that $T$ is absolutely $1$-summing wherever
$T^{\ast}$ almost $p$-summing. So, Corollary \ref{7777} recovers Theorem
\ref{nn11}.

\subsection{Lipschitz $p$-summing and $p$-dominated operators\label{ss2}}

Let $X$ be a pointed metric space with a base point denoted by $0$ and let $Y$ be
an $\mathcal{L}_{p}$-space, with $2\leq p<\infty$. We can easily note that $T$
is Lipschitz $p$-summing if, and only if, it is $RS$-abstract $p$-summing.
Considering $B_{X^{\#}}$ with the pointwise convergence topology,
\[
\left(  K,\mathcal{H}\left(  X;Y\right)  \right)  =\left(  B_{X^{\#}}%
,Lip_{0}(X,Y)\right)  ,
\]%
\[
\left\{
\begin{array}
[c]{c}%
R\colon K\times X\times X\longrightarrow\lbrack0,\infty)\\
R(\varphi,x,q)=|\varphi(x)-\varphi(q)|,
\end{array}
\right.
\]
and
\[
\left\{
\begin{array}
[c]{c}%
S\colon Lip_{0}\left(  X,Y\right)  \times X\times X\longrightarrow
\lbrack0,\infty)\\
S(T,x,q)=\Vert T(x)-T(q)\Vert,
\end{array}
\right.
\]
and considering $X^{d}=X^{\#}$ we can invoke Theorem \ref{q1} to prove the
following extension of the Chen--Zheng Theorem:

\begin{theorem}
Let $X$ be a pointed metric space and $Y$ be an $\mathcal{L}_{p}$-space, with
$2\leq p<\infty$. If $T\in Lip_{0}(X;Y)$ is such that $T^{\#}|_{Y^{\ast}}$ is
almost $p$-summing, then $T$ is Lipschitz $1$-summing.
\end{theorem}

According to Chen and Zheng \cite{ChZh}, a Lipschitz operator $T\colon
X\longrightarrow Y$ between Banach spaces is \textit{Lipschitz }%
$p$\textit{-dominated} if there exist a Banach space $Z$ and an absolutely
$p$-summing linear operator $L\colon X\longrightarrow Z$ such that
\begin{equation}
\left\Vert T(x)-T(q)\right\Vert \leq\left\Vert L(x)-L(q)\right\Vert
\mathrm{~for~all~}x,q\in X. \label{dominated}%
\end{equation}
As a consequence of our main theorem we also conclude that if a Lipschitz
operator $T\colon X\longrightarrow Y$ satisfies (\ref{dominated}) for a
certain $L:X\longrightarrow Z$, where $Z$ is an $\mathcal{L}_{p}$-space, with
$2\leq p<\infty$, then $T$ is Lipschitz $1$-dominated whenever $L^{\ast}$ is
almost $p$-summing.

\subsection{Absolutely $p$-summing $\Sigma$-operators}

According to Angulo-L\'{o}pez and Fern\'{a}ndez-Unzueta \cite{mexicana,
segre}, given Banach spaces $X_{1},\ldots,X_{n}$, the set
\[
\Sigma_{X_{1},\ldots,X_{n}}:=\{x_{1}\otimes\cdots\otimes x_{n}\in X_{1}%
\otimes\cdots\otimes X_{n}:x_{i}\in X_{i},i=1,\ldots,n\}
\]
is the metric space of decomposable tensors endowed with the metric induced by
the projective tensor norm. It is called the \textit{metric Segre cone} of
$X_{1},\ldots,X_{n}$.

If $T: X_{1} \times\cdots\times X_{n} \longrightarrow Y$ is a multilinear
mapping and $Y$ is a vector space, we denote by $\hat{T} \in L(X_{1}
\otimes\cdots\otimes X_{n};Y)$ the unique linear mapping satisfying that for
every $x_{i} \in X_{i}$, with $i \in\{1,\ldots,n\}$ and $T(x_{1},\ldots,x_{n})
= \hat{T}(x_{1}\otimes\cdots\otimes x_{n})$.

\begin{definition}
[see \cite{mexicana}]If $X_{1},\ldots,X_{n}$ are Banach spaces and $Y$ is a
vector space, an operator $f: \Sigma_{X_{1},\ldots,X_{n}} \longrightarrow Y$
is a $\Sigma$-operator if there exists a multilinear operator $T
\in\mathcal{L}(X_{1},..., X_{n},Y)$ such that $f=\hat{T}|_{\Sigma
_{X_{1},\ldots,X_{n}}}$.
\end{definition}

We denote
\[
L(\Sigma_{X_{1},\ldots,X_{n}};Y) = \{f: \Sigma_{X_{1},\ldots,X_{n}}
\longrightarrow Y ~:~ f ~\text{is a}~ \Sigma\text{-operator} \},
\]
and by $\mathcal{L}(\Sigma_{X_{1},\ldots,X_{n}})$ we denote the space of
scalar-valued continuous $\Sigma$-operators endowed with the Lipschitz norm,
which happens to be a dual Banach space in which we shall consider the
weak-star topology.

\begin{definition}
[see \cite{mexicana}]Let $X_{1},\ldots,X_{m},Y$ be Banach spaces. A bounded
$\Sigma$-operator $f\colon\Sigma_{X_{1},\dots,X_{m}}\longrightarrow Y$ is
\textit{absolutely }$p$\textit{-summing} if there is a constant $C$ so that
\[
\left(  \sum\limits_{i=1}^{n}\left\Vert f(x_{i})-f(q_{i})\right\Vert
^{p}\right)  ^{\frac{1}{p}}\leq C\sup_{\varphi\in B_{\mathcal{L}(\Sigma
_{X_{1},\ldots,X_{m}})}}\left(  \sum\limits_{i=1}^{n}\left\vert \varphi
(x_{i})-\varphi(q_{i})\right\vert ^{p}\right)  ^{\frac{1}{p}},
\]
for every natural number $n$ and all $x_{i},q_{i}\in\Sigma_{X_{1},\ldots
,X_{m}}$, with $i=1,...,n$.
\end{definition}

Let us choose, in our abstract framework, $X$ as the pointed metric space
$\Sigma_{X_{1},\ldots,X_{m}}$ with the base point $0=0\otimes\cdots\otimes0$,
and
\[
\left(  K,\mathcal{H}\left(  X;Y\right)  \right)  =\left(  B_{\mathcal{L}%
(\Sigma_{X_{1},\ldots,X_{m}})},Lip_{0}(X,Y)\right)  ,
\]%
\[
\left\{
\begin{array}
[c]{c}%
R\colon K\times X\times X\longrightarrow\lbrack0,\infty)\\
R(\varphi,x,q)=|\varphi(x)-\varphi(q)|,
\end{array}
\right.
\]
and
\[
\left\{
\begin{array}
[c]{c}%
S\colon Lip_{0}\left(  X,Y\right)  \times X\times X\longrightarrow
\lbrack0,\infty)\\
S(f,x,q)=\Vert f(x)-f(q)\Vert.
\end{array}
\right.
\]
Letting $Y$ be an $\mathcal{L}_{p}$-space, with $2\leq p<\infty$,
$X^{d}=X^{\#}=\mathcal{L}(\Sigma_{X_{1},\ldots,X_{m}})$ we have that a bounded
$\Sigma$-operator $f\in Lip_{0}\left(  X;Y\right)  $ is absolutely $p$-summing
if, and only if, it is $RS$-abstract $p$-summing. As $f^{a}=f^{\#}$, Theorem
\ref{q1} and Corollary \ref{7777} provide Kwapie\'{n}-type theorems for this
class of operators.

\subsection{Lipschitz $p$-summing operators at one point}

The next definition is motivated by the notion of absolutely summing arbitrary
mappings between Banach spaces (see \cite{joed} and the references therein).

\begin{definition}
Let $X$ be a normed vector space and $Y$ be a Banach space. A Lipschitz
operator $T:X\longrightarrow Y$ is \textit{Lipschitz} $p$-\textit{summing at
the point} $w\in X$ if there is a constant $C$ such that
\begin{align*}
\left(  \sum_{i=1}^{n}\left\Vert T(w+x_{i})-T(w+q_{i})\right\Vert ^{p}\right)
^{\frac{1}{p}}\leq C\sup_{\varphi\in B_{X^{\#}}}\left(  \sum_{i=1}%
^{n}\left\vert \varphi(x_{i})-\varphi(q_{i})\right\vert ^{p}\right)
^{\frac{1}{p}} , \label{mmu333}%
\end{align*}
for every positive integer $n$ and every $x_{1},...,x_{n},q_{1},...,q_{n}\in
X.$
\end{definition}

Note that Lipschitz $p$-summability at $0$ is precisely the notion of
Lipschitz $p$-summability (see (\ref{765765})). Given $T\in Lip(X,Y)$, we
define $T_{w}:X\longrightarrow Y$ by
\[
T_{w}(x)=T(w+x)-T(w)
\]
and it is simple to check that $T\in Lip(X,Y)$ if, and only if, $T_{w}\in
Lip(X,Y)$ and $Lip(T)=Lip(T_{w})$. Note also that
\[
T\ \text{is Lipschitz}\ p\text{-summing at}\ w\ \text{if, and only if, }%
T_{w}\ \text{is Lipschitz }p\text{-summing}%
\]
and
\[
T\ \text{is Lipschitz}\ p\text{-summing}\ \text{if, and only if, }%
T_{0}\ \text{is Lipschitz }p\text{-summing}.
\]
Considering
\[
w-Lip(X,Y):=\{T_{w}:T\in Lip(X,Y)\},
\]
and choosing $K=B_{X^{\#}}$, with the pointwise convergence topology,
${\mathcal{H}}(X;Y)=w-Lip(X,Y)$ and $R$, $S$ defined by%
\[
\left\{
\begin{array}
[c]{c}%
R\colon B_{X^{\#}}\times X\times X\longrightarrow\lbrack0,\infty)\\
R(\varphi,x,q)=|\varphi(x)-\varphi(q)|,
\end{array}
\right.
\]
and%
\[
\left\{
\begin{array}
[c]{c}%
S\colon w-Lip(X,Y)\times X\times X\longrightarrow\lbrack0,\infty)\\
S(T_{w},x,q)=\left\Vert T_{w}(x)-T_{w}(q)\right\Vert ,
\end{array}
\right.
\]
we conclude that a $T_{w}\in w-Lip(X,Y)$ is $RS$-abstract $p$-summing if, and
only if, $T_{w}$ is Lipschitz $p$-summing (or equivalently, $T$ is Lipschitz
$p$-summing at $w$).

If $Y$ is an $\mathcal{L}_{p}$-space, with $2\leq p<\infty$, $X^{d}=X^{\#}$,
then the abstract adjoint $(T_{w})^{a}$ coincides with the Lipschitz adjoint
$(T_{w})^{\#}$. Consequently, we have the following theorem

\begin{theorem}
Let $X$ be a normed vector space, $Y$ be an $\mathcal{L}_{p}$-space, with
$2\leq p<\infty$, and $T\in Lip(X;Y)$. If $(T_{w})^{\#}|_{Y^{\ast}}$ is almost
$p$-summing, then $T$ is Lipschitz $1$-summing at $w$.
\end{theorem}

Note that the Chen--Zheng theorem generalizes the result of Kwapie\'{n}, since it shows that when $X$ is a normed space we
can replace $Lip_{0}(X;Y)$ by $Lip(X;Y).$


\begin{thebibliography}{99}                                                                                               %


\bibitem {mexicana}J.C. Angulo-L\'opez and M. Fern\'andez-Unzueta,
\textit{Lipschitz $p$-summing multilinear operators}, to appear in J. Funct. Anal.

\bibitem {bayart}F. Bayart, \textit{Multiple summing maps: coordinatewise
summability, inclusion theorems and $p$-Sidon sets}, J. Funct. Anal.,
\textbf{274} (2018), 1129--1154.

\bibitem {BBJ}G. Botelho, H.A. Braunss and H. Junek, \textit{Almost
$p$-summing polynomials and multilinear mappings}, Arch. Math, \textbf{76}
(2001), 109--118.

\bibitem {jfa}G. Botelho, D. Pellegrino and P. Rueda, \textit{Pietsch's
factorization theorem for dominated polynomials}, J. Funct. Anal.,
\textbf{243} (2007), 257--269.

\bibitem {JA}J.A. Ch\'{a}vez-Dom\'{\i}nguez, \textit{Duality for Lipschitz
p-summing operators}, J. Funct. Anal., \textbf{261} (2011), 387--407.

\bibitem {JAproc}J.A. Ch\'{a}vez-Dom\'{\i}nguez, \textit{Lipschitz $(q,
p)$-mixing operators}, Proc. Amer. Math. Soc., \textbf{140} (2012), 3101--3115.

\bibitem {CZ}D. Chen and B. Zheng, \textit{Remarks on Lipschitz $p$-summing
operators}, Proc. Amer. Math. Soc., \textbf{139} (2011), 2891--2898.

\bibitem {ChZh}D. Chen and B. Zheng, \textit{Lipschitz $p$-integral operators
and Lipschitz $p$-nuclear operators}, Nonlinear Anal., \textbf{75} (2012), 5270--5282.

\bibitem {djt}J. Diestel, H. Jarchow and A. Tonge, \textit{Absolutely Summing
Operators}, Cambridge Studies in Advanced Mathematics, Cambridge University
Press, 1995.

\bibitem {FJ}J. Farmer and W. B. Johnson, \textit{Lipschitz $p$-summing
operators}, Proc. Amer. Math. Soc., \textbf{137} (2009), 2989--2995.

\bibitem {Fer}M. Fern\'{a}ndez-Unzueta, \textit{Dunford-Pettis and
Dieudonn\'{e} polynomials on Banach spaces}, Illinois J. Math., \textbf{45}
(2001), 291--307.

\bibitem {segre}M. Fern\'{a}ndez-Unzueta, \textit{The Segre cone of Banach
spaces and multilinear mappings}, Linear Multilinear Algebra, \textbf{68}
(2020), 575--593.

\bibitem {gro}A. Grothendieck, \textit{R\'{e}sum\'{e} de la th\'{e}orie
m\'{e}trique des produits tensoriels topologiques}. (French) Bol. Soc. Mat.
S\~ao Paulo, \textbf{8} (1953) 1--79.

\bibitem {K}S. Kwapie\'{n}, \textit{On a theorem of L. Schwartz and the
applications to absolutely summing operators}, Studia Math., \textbf{38}
(1970), 193--201.

\bibitem {LP}J. Lindenstrauss, A. Pe{\l }czy{\'n}ski, \textit{Absolutely
summing operators in $L_{p}$-spaces and their applications}, Studia Math.,
\textbf{29} (1968), 275--326.

\bibitem {MM}L. Maligranda and M. Masty\l {}o, \textit{Inclusion mappings
between Orlicz sequence spaces}, J. Funct. Anal., \textbf{176} (2000), 264--279.

\bibitem {joed}D. Pellegrino and J. Santos, \textit{On summability of
nonlinear maps: a new approach}, Math. Z., \textbf{270} (2012), 189--196.

\bibitem {adv}D. Pellegrino, J. Santos and J. B. Seoane-Sep\'{u}lveda,
\textit{Some techniques on nonlinear analysis and applications}, Adv. Math.,
\textbf{229} (2012), 1235--1265.

\bibitem {LONDON}D. Pellegrino, J. Santos and J. B. Seoane-Sep\'{u}lveda,
\textit{A general extrapolation theorem for absolutely summing operators},
Bull. Lond. Math. Soc., \textbf{44} (2012), 1292--1302.

\bibitem {PSDianaE}D. Pellegrino, J. Santos, D. Serrano-Rodr\'{\i}guez and E.
Teixeira, \textit{A regularity principle in sequence spaces and applications},
Bull. Sci. Math., \textbf{141} (2017), 802--837.

\bibitem {perez}D. P\'{e}rez-Garc\'{\i}a, \textit{Comparing different classes
of absolutely summing multilinear operators}, Arch. Math., \textbf{85} (2005), 258--267.

\bibitem {Pie}A. Pietsch, \textit{Absolut $p$-summierende Abbildungen in
normieten R\"{a}umen}, Studia Math., \textbf{27} (1967), 333--353.

\bibitem {PPPP}A. Pietsch. \textit{Ideals of multilinear functionals},
Proceedings of the Second International Conference on Operator Algebras,
Ideals and Their Applications in Theoretical Physics, 185-199, Teubner-Texte,
Leipzig, 1983.

\bibitem {Sawashima}I. Sawashima, \textit{Methods of duals in nonlinear
analysis}, Lecture Notes in Econom. and Math. Systems, \textbf{419} (1995), 247--259.

\bibitem {vieira}F. Vieira Costa Junior, \textit{The optimal multilinear
Bohnenblust-Hille constants: a computational solution for the real case},
Numer. Funct. Anal. Optim., \textbf{39} (2018), 1656--1668.
\end{thebibliography}
\end{document}